\documentclass[12pt]{article}
\usepackage{amssymb}

\newcounter{enucount}
\newcounter{enuref}
\newcounter{enualph}
\newcounter{enunumb}
\renewcommand{\theenucount}{(\roman{enucount})}

\renewcommand{\theenualph}{(\alph{enualph})}

\newtheorem{sdef}{Definition}[section]

\newtheorem{sthm}[sdef]{Theorem}
\newtheorem{slem}[sdef]{Lemma}

\newtheorem{scor}[sdef]{Corollary}
\newtheorem{sques}[sdef]{Question}

\newtheorem{sclaim}[sdef]{Claim}
\newtheorem{smainclaim}[sdef]{Main Claim}

\newenvironment{otherthm}[1]{\bigskip \noindent {\bf #1} \em}
   {\bigskip}
\newenvironment{proof}{\noindent\emph{Proof.}\enspace}%
   {\hfill$\square$\smallskip}

\newenvironment{alphenumerate}{\begin{list}{\rm \theenualph}{\usecounter{enualph}
    \setlength{\labelwidth}{1cm}}}
   {\end{list}}

\newcommand{\forces}{\Vdash}
\newcommand{\re}{{\upharpoonright}}

\newcommand{\F}{{\cal F}}

\newcommand{\PP}{{\mathbb P}}
\newcommand{\QQ}{{\mathbb Q}}

\renewcommand{\SS}{{\mathbb S}}

\newcommand{\bb}{{\mathfrak b}}
\newcommand{\cc}{{\mathfrak c}}

\newcommand{\ee}{{\mathfrak e}}

\newcommand{\mm}{{\mathfrak m}}

\newcommand{\sss}{{\mathfrak s}}
\newcommand{\vv}{{\mathfrak v}}

\newcommand{\Add}[1]{{{\mathsf{add}}({\cal #1})}}

\newcommand{\Cof}[1]{{{\mathsf{cof}}({\cal #1})}}

\newcommand{\stem}{{\mathrm{stem}}}

\newcommand{\dom}{{\mathrm{dom}}}

\newcommand{\const}{{\mathrm{const}}}
\newcommand{\linked}{{\mathrm{linked}}}

\newcommand{\sub}{\subseteq}
\newcommand{\sem}{\setminus}
\newcommand{\twoom}{2^\omega}
\newcommand{\twolom}{2^{<\omega}}

\newcommand{\ha}{\,{}\hat{}\,}
\newcommand{\la}{\langle}
\newcommand{\ra}{\rangle}

\newcommand{\noi}{\noindent}

\begin{document}

\title{Evasion and prediction IV \\
Fragments of constant prediction}

\author{J\"org Brendle\thanks{Supported by Grant--in--Aid for Scientific Research
   (C)(2)12640124, Japan Society for the Promotion of Science}  \\ 
   The Graduate School of Science and Technology \\
   Kobe University \\
   Rokko--dai 1--1, Nada--ku \\
   Kobe 657--8501, Japan 
\and
   Saharon Shelah  \\
   Institute of Mathematics \\
   The Hebrew University of Jerusalem \\
   91904 Jerusalem, Israel \\
   and \\
   Department of Mathematics \\
   Rutgers University \\
   New Brunswick, NJ 08903, USA}

\maketitle

\begin{abstract}
\noi Say that a function $\pi : n^{< \omega} \to n$ (henceforth
called a predictor) $k$--constantly predicts a real $x \in  n^\omega$ if for almost all
intervals $I$ of length $k$, there is $i \in I$  such that $x(i) = \pi (x \re i)$.
We study the $k$--constant prediction number $\vv_n^\const (k)$, that is, the size
of the least family of predictors needed to $k$--constantly predict all reals, for different
values of $n$ and $k$, and investigate their relationship. 
\end{abstract}


\section*{Introduction}

This work is about evasion and prediction, a combinatorial concept originally introduced
by Blass when studying set--theoretic aspects of the Specker phenomenon in abelian
group theory~\cite{Bl94}. The motivation for our investigation came from a (still open)
question of Kamo, as well as from an argument in a proof by the first author. Let us explain this
in some detail.

For our purposes, let $n\leq\omega$ and call a function $\pi : n^{< \omega} \to n$ a {\it predictor}.
Say $\pi$ {\it $k$--constantly predicts} a real $x \in  n^\omega$ if for almost all
intervals $I$ of length $k$, there is $i \in I$  such that $x(i) = \pi (x \re i)$.
In case $\pi$ $k$--constantly predicts $x$ for some $k$, say that $\pi$ {\it constantly
predicts} $x$. The {\it constant prediction number $\vv^\const_n$}, introduced by Kamo in
\cite{Ka00}, is the smallest
size of a set of predictors $\Pi$ such that every $x \in n^\omega$ is constantly predicted by some
$\pi\in\Pi$. Kamo~\cite{Ka00} showed that $\vv_\omega^\const$ may be larger than
all the $\vv_n^\const$ where $n \in \omega$. He asked

\begin{otherthm}{Question.} {\rm (Kamo~\cite{Kata})}
Is $\vv_2^\const = \vv_n^\const$ for all $n\in\omega$.
\end{otherthm}

Some time ago, the first author answered another question of Kamo's by showing that $\bb \leq \vv_2^\const$
where $\bb$ is the unbounding number~\cite{Brta}. Now, the standard approach to such a result would have been to show that,
given a model $M$ of $ZFC$ such that there is a dominating real $f$ over $M$, there must be a real which
is not constantly predicted by any predictor from $M$. This, however, is far from being true.
In fact, one needs a sequence of $2^k -1$ models $M_i$ and dominating reals $f_i$ over $M_i$
belonging to $M_{i+1}$ to be able to construct a real which is
not $k$--constantly predicted by any predictor from $M_0$, and this result is optimal (see~\cite{Brta}
for details). This means $k$--constant prediction gets easier in a strong sense the larger
$k$ gets, and one can expect interesting results when investigating the cardinal invariants which
can be distilled out of this phenomenon.

Accordingly, let us define the {\it $k$--constant prediction number $\vv^\const_n (k)$} to be the size of
the smallest set of predictors $\Pi$ such that every $x \in n^\omega$ is $k$--constantly predicted
by some $\pi\in\Pi$.
Interestingly enough, Kamo's question cited above has a positive answer when relativized to the
new situation. Namely, we shall show in Section~\ref{1sec} that $\vv_2^\const (k) = \vv_n^\const (k)$
for all $k, n < \omega$ (see~\ref{1mainthm}). 
Moreover, for $k < \ell$, one may well have $\vv_2^\const (\ell) < \vv_2^\const (k)$ (Theorem~\ref{Sacksprod}).
Any hope to use Theorem~\ref{1mainthm} as an intermediate step to answer Kamo's question is dashed,
however, by Theorem~\ref{2mainthm} which says that $\vv_2^\const$ may be strictly smaller than
the minimum of all $\vv_2^\const (k)$'s.

In Section~\ref{3sec}, we dualize Theorem~\ref{Sacksprod} to a consistency result about evasion numbers
and establish a connection between those and Martin's axiom for $\sigma -k$--linked partial
orders (see Theorem~\ref{3mainthm}). 

We keep our notation fairly standard. For basics concerning the cardinal invariants considered here,
as well as the forcing techniques, see~\cite{BJ95} and~\cite{Blta}.

The results in this paper were obtained in September 2000 during and shortly after the
second author's visit to Kobe. The results in Sections~\ref{1sec} and \ref{2sec}
are due to the second author. The remainder is the first author's work.


\section{The $ZFC$--results}
\label{1sec}

Temporarily say that $\pi : n^{<\omega} \to n$ {\it weakly $k$--constantly predicts}
$x \in n^\omega$  if for almost all 
$m$ there is $i < k$ such that $\pi (x \re mk+ i)
= x (mk + i)$. This notion is obviously weaker than $k$--constant prediction.
It is often more convenient, however.
We shall see soon that in terms of cardinal invariants the two notions are the same.

Put $G = \{ \bar g  = \la g_i ; \; i < k  \ra ; \; g_i : n^k \to 2 \}$.

\begin{sthm}   \label{sec1mainthm}
There are functions $\bar \pi = \la \pi^{\bar g, j} ; \; (\bar g,j) \in G \times k \ra \mapsto \psi_{\bar\pi}$ (where 
$\pi^{\bar g,j} : \twolom \to 2$ and $\psi_{\bar\pi} : n^{<\omega} \to n$)
and $y  \mapsto \la y^{\bar g , j} ; \; (\bar g,j) \in G \times k \ra$ 
(where $y\in n^\omega$ and $y^{\bar g,j} \in \twoom$)
such that if $\pi^{\bar g , j}$ weakly $k$--constantly predicts $y^{\bar g,j}$  for
all pairs $(\bar g, j)$,
then $\psi_{\bar\pi}$ $k$--constantly predicts $y$.
\end{sthm}

\begin{proof}
Given $y \in n^\omega$, define $y^{\bar g,j}$ by
\[ y^{\bar g,j} (mk  + i) = g_i (y \re [mk + j , (m+1) k + j) ). \]
Also, for $\sigma \in n^{< \omega}$, say $|\sigma| = m_0 k + j$, define $\sigma^{\bar g,j}$ by
\[ \sigma^{\bar g, j} (mk + i) =  g_i (\sigma \re [mk + j , (m+1) k + j) ) \]
for all $m < m_0$. So $|\sigma^{\bar g,j} | = m_0 k$.

Given $\bar \pi =  \la \pi^{\bar g, j} ; \; (\bar g,j) \in G \times k \ra$, a sequence of predictors
for the space $\twoom$, and $\sigma \in n^{< \omega}$, say $|\sigma| = m k + j$,
put 
\[A_\sigma^k = \{ \tau \supset \sigma; \; |\tau| = |\sigma| + k 
\mbox{ and } \forall \bar g \; \exists i \; (\tau^{\bar g,j} (mk  + i) 
= \pi^{\bar g,j} (\tau^{\bar g,j} \re mk+i) ) \}. \]
For $i <k$, define $A_\sigma^i = \{ \tau \supset \sigma ; \; \tau \in A_{\sigma \re |\sigma| - k + i}^k \}$.
So, if $\tau \in A_\sigma^i$, $|\tau| = |\sigma| + i$.


\begin{sclaim}
$|A_\sigma^k| < 2^k$ for all $\sigma$.
\end{sclaim}

\begin{proof}
Assume that, for some $\sigma$, we have $|A_\sigma^k| \geq 2^k$.
List $\{ \tau_\ell ; \; \ell < 2^k \} \sub A_\sigma^k$ and list
$2^k = \{ \sigma_\ell ; \; \ell < 2^k \}$.
Fix $m$ and $j$ such that $|\sigma| = mk + j$. 
Define $g_i (\tau_\ell \re [mk+j, (m+1) k + j)) = \sigma_\ell (i)$ and consider $\bar g = \la g_i ; \; i < k \ra$.
Then $\tau_\ell^{\bar g, j} \re [ mk  , (m+1) k ) = \sigma_\ell$.
This is a contradiction to the definition of $A_\sigma^k$ for it would mean $\pi^{\bar g,j}$ cannot predict correctly
all $\tau_\ell^{\bar g, j}$ somewhere in the interval $[mk  , (m+1) k )$.
\end{proof}

\medskip

For $\sigma \in n^{<\omega}$ define $\psi_{\bar\pi} (\sigma)$ as follows.
First let $i \leq k$ be minimal such that $|A_\sigma^i| < 2^i$. Such $i$ exists by the claim.
Then let $\psi_{\bar\pi} (\sigma)$ be any $\ell$ such that $A_{\sigma\ha\la\ell\ra}^{i-1}$ is of maximal
size.

To see that this works, let $y \in n^\omega$. Let  $\pi^{\bar g,j}$ be predictors such that
for all $\bar g , j$ and almost all $m$, there is $i$ such that $y^{\bar g,j} (mk  + i) = \pi^{\bar g,j} 
(y^{\bar g,j} \re mk+i) $. Fix $m_0$ such that for all $m \geq m_0$ and all $\bar g,j$,
there is $i$ such that $y^{\bar g,j} (mk  + i) = \pi^{\bar g,j} 
(y^{\bar g,j} \re mk+i) $.
Let $mk+j \in\omega$ with $m \geq m_0$. Thus $y \re mk +j+i \in A^i_{y \re mk+j}$
for all $i \leq k$. We need to find $i<k$ such that $\psi_{\bar\pi} (y\re mk+j+i) = y(mk+j+i)$.
To this end simply note that if $i$ is such that $\psi_{\bar\pi} (y\re mk+j+i) \neq y(mk+j+i)$,
then, by definition of $\psi_{\bar\pi}$,
\[ |A_{y\re mk+j+i+1}^{\ell_i -1} | \leq {|  A_{y\re mk+j+i}^{\ell_i}   | \over 2}\] where $\ell_i$ is minimal
with $|  A_{y\re mk+j+i}^{\ell_i}   | < 2^{\ell_i}$. This means in particular $|A_{y\re mk+j+i+1}^{\ell_i -1} |
< 2^{\ell_i -1}$. A fortiori, $\ell_{i+1} \leq \ell_i - 1$. Since $\ell_0 \leq k$,
this entails that if we had $\psi_{\bar\pi} (y\re mk+j+i) \neq y(mk+j+i)$ for all $i<k$, we would
get $\ell_i = 0$ for some $i \leq k$. Thus  $|  A_{y\re mk+j+i}^{0}   | < 2^{0} =1$.
So $A_{y\re mk+j+i}^{0} = \emptyset$. However $y\re mk+j+i \in A_{y\re mk+j+i}^{0}$, a contradiction.
This completes the proof of the theorem.
\end{proof}

\medskip

Define the {\it $k$--constant evasion number $\ee^{\const}_n (k)$} to be the dual of $\vv^{\const}_n (k)$,
namely the size of the smallest set of functions $F \sub n^\omega$ such that for every predictor
$\pi$ there is $x \in F$ which is no $k$--constantly predicted by $\pi$.
Similarly, define the {\it constant evasion number} $\ee^{\const}_n $.

Let $\bar \vv_n^\const (k)$ denote the size of the least family $\Pi$ of predictors
$\pi: n^{<\omega} \to n$ such that every $y \in n^\omega$ is weakly $k$--constantly
predicted by a member of $\Pi$. Dually, $\bar \ee_n^\const (k)$ is the size of the least
family $F \sub n^\omega$ such that no predictor $\pi : n^{<\omega} \to n$
weakly $k$--constantly predicts all members of $F$. The above theorem entails

\begin{scor}
$\vv_n^\const (k) \leq \bar \vv_2^\const (k)$. Dually, $\ee_n^\const (k) \geq \bar \ee_2^\const (k)$.
\end{scor}

\begin{proof}
Let $\Pi$ be a family of predictors in $\twoom$ weakly $k$--constantly predicting all functions.
Put $\Psi = \{ \psi_{\bar\pi} ; \; \bar\pi = \la \pi^{\bar g,j} ; \; (\bar g,j) \in G\times k \ra \in \Pi^{<\omega} \}$.
By the theorem, every $y \in n^\omega$ is $k$--constantly predicted by a member of $\Psi$. 
This shows $\vv_n^\const (k) \leq \bar \vv_2^\const (k)$.

Next let $F \sub n^\omega$ be a family of functions such that no predictor
$k$--constantly predicts all of $F$.
Let $Y = \{ y^{\bar g, j} ; \; (\bar g,j) \in G\times k$ and $y\in F \} \sub 2^\omega$.
Assume $\pi : \twolom \to 2$ weakly $k$--constantly predicts all members of $Y$.
Then $\psi_{\bar\pi} $ $k$--constantly predicts all members of $F$,
where we put $\bar\pi = \la  \pi^{\bar g,j} ; \; (\bar g,j) \in G\times k \ra$ with $\pi^{\bar g,j}
= \pi$ for all $(\bar g ,j) \in G \times k$, a contradiction.
\end{proof}

\medskip

Since the other inequalities are trivial, we get

\begin{sthm} \label{1mainthm}
$\bar \vv_n^\const (k) = \vv_n^\const (k) = \vv_2^\const (k)$ for all $n$.
Dually, $\bar \ee_n^\const (k) = \ee_n^\const (k) = \ee_2^\const (k)$ for all $n$.
\end{sthm}

A fortiori, we also get $\min\{ \vv_n^\const (k) ; \; k \in \omega \} = \min \{ \vv_2^\const (k) ; \; k \in \omega \}$
and $\sup \{ \ee_n^\const (k) ; \; k \in \omega \} = \sup \{ \ee_2^\const (k) ; \; k \in \omega \}$ for all $n$.


\section{Prediction and relatives of Sacks forcing}
\label{2sec}

For $2 \leq k < \omega$, define
{$k$--ary Sacks forcing $\SS^k$} to be the set of all subtrees $T \sub k^{<\omega}$
such that below each node $s \in T$, there is $t \supset s$ whose $k$ immediate
successor nodes $t\ha\la i \ra$ ($i<k$) all belong to $T$. $\SS^k$ is ordered by inclusion.
Obviously $\SS^2$ is nothing but standard Sacks forcing $\SS$.

Iterating  $\SS^k$ $\omega_2$ many times with countable support over
a model for $CH$ yields a model where $\vv_2^\const (\ell)$ is large if $2^\ell \leq k$
and small otherwise. This has been observed independently around the same time
by Kada~\cite{Kadta}. However, one can get better consistency results 
by using large countable support products instead. The following is in the
spirit of~\cite{GoSh93}.

\begin{sthm}  \label{Sacksprod}
Assume $CH$. Let $2 \leq k_1 < ... < k_{n-1}$. Also let $\kappa_i$, $i \leq n$,
be cardinals with $\kappa_i^\omega = \kappa_i$ and $\kappa_{n} < ... < \kappa_0$.
Then there is a generic extension satisfying $\vv^\const_2 = \min \{ \vv^\const_2 (k) ; \; k\in\omega\} = \vv^\const_2
(k_{n-1} +1) = \kappa_{n}  $, $\vv_2^\const (k_i) = \vv^\const_2 (k_{i-1} +1) = \kappa_{i}  $ for $0< i < n$
and $\cc = \kappa_0$.
\end{sthm}

\begin{proof}
We force with the countable support product $\PP = \prod_{\alpha < \kappa_0} \QQ_\alpha$ where 
\begin{itemize}
\item $\QQ_\alpha$ is Sacks forcing $ \SS_\alpha$ for $\kappa_1\leq\alpha < \kappa_0$,
\item $\QQ_\alpha$ is $2^{k_i}$--ary Sacks forcing $\SS_\alpha^{2^{k_i}}$ for $0<i<n$ and 
   $\kappa_{i+1}\leq\alpha < \kappa_i$, and
\item $\QQ_\alpha$ is $\SS_\alpha^{\ell_\alpha}$ where $|\{ \alpha ; \; \ell = \ell_\alpha \} | = \kappa_n$
   for all $\ell$, for $\alpha < \kappa_n$.
\end{itemize}
By $CH$, $\PP$ preserves cardinals and cofinalities.
$\cc = \kappa_0$ is also immediate.

Note that if $X \sub \twoom$ and $|X| < \kappa_i$, then there is $A \sub \kappa_0$ of size $< \kappa_i$
such that $X \in V [G_A]$, the generic extension by conditions with support contained in $A$,
i.e. via the ordering $\prod_{\alpha \in A} \QQ_\alpha$.
So there is $\alpha \in (\kappa_i \sem \kappa_{i+1}) \sem A$. 
Clearly the generic real added by $\QQ_\alpha = \SS_\alpha^{2^{k_i}}$  is not $k_i$--constantly
predicted by any predictor from $V[G_A]$. This shows $\vv_2^\const (k_i) \geq \kappa_i$.
A similar argument shows $\vv_2^\const \geq \kappa_n$.

So it remains to see that $\vv_2^\const (k_{i_0 -1} +1) \leq \kappa_{i_0}$ for $0 < i_0 \leq n$.
Put $\ell = k_{i_0 -1} +1$. Let $\dot f$ be a $\PP$--name for a function in $2^\omega$. 
By a standard fusion argument we can recursively construct 
\begin{itemize}
\item a strictly increasing sequence $m_j$, $j\in\omega$,
\item $A \sub \kappa_0$ countable,
\item $\la D_\alpha ; \; \alpha \in A \ra$, a partition of $\omega$ into countable sets,
\item a condition $p = \la p_\alpha ; \; \alpha \in A \ra \in \PP$, and
\item a tree $T \sub \twolom$
\end{itemize}
such that
\begin{alphenumerate}
\item if $\sigma \in T \cap 2^{m_j}$, $j \in D_\alpha$, and $\alpha \in \kappa_i \sem \kappa_{i+1}$ ($i<n$),
   then $| \{ \tau \in T \cap 2^{m_{j+1}} ; \; \sigma \sub \tau \} | = 2^{k_i}$ where we put $k_0 = 1$, \label{prod-side}
\item $p \forces \dot f \in [T]$, and  
\item whenever $q \leq p$ where $q = \la q_\beta ; \; \beta \in B \ra$ with $A \sub B$, 
   $\sigma \in T \cap 2^{m_j}$, and $j \in D_\alpha$ are such that $q \forces \sigma \sub \dot f$,
   then there are $r_\alpha \leq q_\alpha$ and $\tau \in T \cap 2^{m_{j+1}}$ with $\tau \supseteq \sigma$,
   such that $r \forces \tau \sub \dot f$ where $r = \la r_\beta ; \; \beta \in B \ra$ with $r_\beta = q_\beta$
   for $\beta \neq \alpha$.   \label{prod-main}
\end{alphenumerate}
Now let $G_{\kappa_{i_0}}$ be $\prod_{\alpha < \kappa_{i_0}} \QQ_\alpha$--generic with $p\re\kappa_{i_0} 
\in G_{\kappa_{i_0}}$.
By \ref{prod-main} above, there is, in $V[G_{\kappa_{i_0}}]$, a tree $S \sub T$
such that for all $\alpha \in A \cap \kappa_{i_0}$, $j \in D_\alpha$ and $\sigma \in S \cap
2^{m_j}$, there is a unique $\tau \in S \cap 2^{m_{j+1}}$ extending $\sigma$,
and such that $\dot f$ is forced to be a branch of $S$ by the remainder of the forcing below $p$.
By \ref{prod-side}, we also have that for all $\alpha \in A \sem \kappa_{i_0}$, $j \in D_\alpha$ and $\sigma \in S \cap
2^{m_j}$, there are at most $2^{k_{i_0 -1}}$ many $\tau \in S \cap 2^{m_{j+1}}$ extending $\sigma$.
This means we can recursively construct a predictor $\pi \in V[G_{\kappa_{i_0}}]$ which 
$\ell$--constantly predicts all branches of $S$. A fortiori,
$\dot f$ is forced to be predicted by $\pi$ by the remainder of the forcing below $p$.
On the other hand, $V[G_{\kappa_{i_0}}]$ satisfies $\cc = \kappa_{i_0}$ so that there are a total number
of $\kappa_{i_0}$ many predictors in $V[G_{\kappa_{i_0}}]$, and they $\ell$--constantly predict all
reals of the final extension. This completes the argument.
\end{proof}

\medskip

It is easy to see that in models obtained by such product constructions, 
$\vv^\const_2 = \min \{ \vv^\const_2 (k) ; \; k\in\omega\}$ must always hold.
To distinguish between these two cardinals, we must turn once again to a countable support
iteration. 

\begin{sthm}  \label{2mainthm}
Assume $CH$. There is a generic extension satisfying $\vv_2^\const = \aleph_1 < \min \{ \vv_2^\const (k) ;
\; k\in\omega\} = \cc = \aleph_2 $.
\end{sthm}

\begin{proof}
Let $\la k_\alpha  ; \; \alpha < \omega_2 \ra$ be a sequence of natural numbers $\geq 2$ in which each $k$ appears
$\omega_2$ often and such that in each limit ordinal, the set of $\alpha$ with $k_\alpha = 2$
is cofinal. 

We perform a countable support iteration $\la \PP_\alpha, \dot \QQ_\alpha ; \; \alpha < \omega_2
\ra$ such that 
\[ \forces_\alpha ``\dot \QQ_\alpha = \dot \SS^{k_\alpha}, \mbox{ that is } k_\alpha \mbox{--ary Sacks forcing}." \]
By $CH$, $\PP_{\omega_2}$ preserves cardinals and cofinalities. As in the previous proof,
we see $\vv_2^\const (k) = \cc = \aleph_2$ for all $k$. We are left with showing that
$\vv_2^\const = \aleph_1$.

Let $\dot f$ be a $\PP_{\omega_2}$--name for a function in $\twoom$.
Notice given any $p_0 \in \PP_{\omega_2}$, we can find $p \leq p_0$ and $\alpha < \omega_2$
such that
\[ p \forces \dot f \in V[\dot G_\alpha] \sem \bigcup_{\beta<\alpha} V[\dot G_\beta]. \]
First consider the case $\alpha$ is a successor ordinal, say $\alpha = \beta + 1$.
Let $\ell$ be such that $2^\ell > k_\beta$. The following is the main point.

\begin{smainclaim}
There are $q \leq p$ and a predictor $\pi \in V$ such that
\[ q \forces `` \pi \;\; \ell \mbox{--constantly predicts } \dot f."\]
\end{smainclaim}

\begin{proof}
We construct recursively
\begin{itemize}
\item $A \sub \alpha$ countable,
\item $\la D_\gamma ; \; \gamma \in A \ra$, a partition of $\omega$ into countable sets,
\item finite partial functions $a_j : A \to \omega$, $j\in\omega$, 
\item conditions $p_j \in \PP_\alpha$, $j\in\omega$, 
\item a strictly increasing sequence $m_j$, $j\in\omega$, 
\item a tree $T \sub \twolom$, and
\item a predictor $\pi: \twolom \to 2$
\end{itemize}
such that
\begin{alphenumerate}
\item $\beta \in A$,
\item $a_0 = \emptyset$,
\item if $j \in D_\gamma$, then $\dom (a_{j+1}) = \dom (a_j) \cup \{ \gamma \}$; in case $\gamma \notin \dom (a_j)$,
   we have $a_{j+1} (\gamma) =0$, otherwise $a_{j+1} (\gamma) = a_j (\gamma) + 1$;
   $a_{j+1} (\delta) = a_j (\delta)$ for $\delta \neq \gamma$,   \label{iter-main-3}
\item $p_0 = p$,   \label{iter-main-4}
\item $p_{j+1} \leq p_j$; furthermore for all $\gamma \in \dom (a_{j+1})$, \\
   $p_{j+1} \re \gamma \forces_\gamma p_{j+1} (\gamma) \leq_{a_{j+1} (\gamma)} p_j (\gamma)$,  \label{iter-main-5}
\item $\bigcup_j \dom (p_j) = \bigcup_j \dom (a_j) = A$,  \label{iter-main-6}
\item if $\sigma \in T \cap 2^{m_j}$, $j \in D_\gamma$, 
   then $| \{ \tau \in T \cap 2^{m_{j+1}} ; \; \sigma \sub \tau \} | = k_\gamma$, 
\item for each $\sigma \in  T \cap 2^{m_j}$, there is $p^\sigma_j \leq p_j$ which forces $\sigma \sub \dot f$;
   furthermore $p_j \forces \dot f \re m_j \in T \cap 2^{m_j}$, and  \label{iter-main-8}
\item $\pi$ $\ell$--constantly predicts all branches of $T$.   \label{iter-main-9}
\end{alphenumerate}
Most of this is standard. There is, however, one trick involved, and we describe the construction.
For $j = 0$, there is nothing to do. So assume we arrived at stage $j$, and we are supposed to
produce the required objects for $j+1$. This proceeds by recursion on $\sigma \in T \cap 2^{m_j}$.
Since the recursion is straightforward, we confine ourselves to describing a single step.

Fix $\sigma \in T \cap 2^{m_j}$. Let $\gamma$ be such that $j \in D_\gamma$.
Without loss $\gamma < \beta$ (the case $\gamma = \beta$ being easier).
Consider $p_j^\sigma$. Step momentarily into $V[G_\beta]$
with $p_j^\sigma \re \beta \in G_\beta$. Then $p_j^\sigma (\beta) \forces_{ \QQ_\beta}
\sigma \sub \dot f$. Since $\dot f$ is forced not to be in $V[G_\beta]$,
we can find $m^\sigma \in \omega$, pairwise incompatible $r^\sigma_i \leq p_j^\sigma (\beta)$,
and distinct $\tau^\sigma_i \in 2^{m^\sigma}$ where $i < k_\gamma$ extending $ \sigma$ such that
$r^\sigma_i \forces_{\QQ_\beta} \tau^\sigma_i \sub \dot f$. As $\QQ_\beta$
is $k_\beta$--ary Sacks forcing, we may do this in such a way that the predictor
$\pi$ can be extended to $\ell$--constantly predict all $\tau^\sigma_i$.

Back in $V$, by extending the condition $p_j^\sigma$ if necessary, we may without
loss assume that it decides $m^\sigma$ and the $\tau^\sigma_i$.
We therefore have the extension of $\pi$ which $\ell$--constantly predicts all $\tau^\sigma_i$
already in the ground model $V$. We may also suppose that $p_j^\sigma \re \gamma$ decides the stem
of $p_j^\sigma (\gamma)$, say $p_j^\sigma \re \gamma \forces_\gamma \stem (p_j^\sigma (\gamma)) = t$.
For $i < k_\gamma$ define $p_{j+1}^{\tau^\sigma_i}$ such that 
\begin{itemize} 
\item $p_{j+1}^{\tau^\sigma_i} \re \gamma = p_j^\sigma \re \gamma$, $p_{j+1}^{\tau^\sigma_i} \re [\gamma +1,
   \beta ) = p_j^\sigma \re [\gamma +1, \beta )$, 
\item $p_{j+1}^{\tau^\sigma_i} \re \gamma \forces_\gamma p_{j+1}^{\tau^\sigma_i} (\gamma) = (p_j^\sigma (\gamma))_{t \ha
   \la i \ra }$,
\item $p_{j+1}^{\tau^\sigma_i} \re \beta \forces_\beta  p_{j+1}^{\tau^\sigma_i} (\beta) = \dot r^\sigma_i$.
\end{itemize}
Doing this (in a recursive construction)
for all $\sigma \in T \cap 2^{m_j}$ and increasing $m^\sigma$ if necessary, we may assume there
is $m_{j+1}$ with $m_{j+1} = m^\sigma$ for all $\sigma$. Finally $p_{j+1}$ is the least upper bound of 
all the $p_{j+1}^{\tau^\sigma_i}$.

This completes the construction. By \ref{iter-main-3}, \ref{iter-main-5}, and \ref{iter-main-6},
the sequence of $p_j$'s has a lower bound $q \in \PP_\alpha$. By \ref{iter-main-4}, $q \leq p$.
By \ref{iter-main-8}, $q \forces \dot f \in [T]$ which means that \ref{iter-main-9} entails
$q \forces ``\dot f$ is $\ell$--constantly predicted by $\pi$," as required.
\end{proof}

\medskip

Now let $\alpha$ be a limit ordinal.
Using a similar argument and the fact that below $\alpha$, $\dot\QQ_\beta$ is cofinally often Sacks forcing,
we see

\begin{sclaim}
There are $q \leq p$ and a predictor $\pi \in V$ such that
\[ q \forces `` \pi \;\; 2 \mbox{--constantly predicts } \dot f."\]
\end{sclaim}

This completes the proof of the theorem.
\end{proof}


\section{Evasion and fragments of $MA(\sigma$--linked)}
\label{3sec}

Let $k \geq 2$.
Recall that a partial order $\PP$ is said to be {\it $\sigma -k$--linked} if it can be written as a countable
union of sets $P_n$ such that each $P_n$ is $k$--linked, that is, any $k$ many elements from
$P_n$ have a common extension. Clearly every $\sigma$--centered forcing is $\sigma -k$--linked for
all $k$, and a $\sigma -k$--linked p.o. is also $\sigma - (k-1)$--linked.
Random forcing is an example of a p.o. which is $\sigma - k$--linked for all $k$, yet not
$\sigma$--centered. A p.o. with the former property shall be called
{\it $\sigma - \infty$--linked} henceforth. We shall deal  with p.o.'s which arise naturally in connection with
constant prediction and which are $\sigma - (k-1)$--linked but not $\sigma -k$--linked for some $k$.
Let $\mm (\sigma-k$--linked) denote the least cardinal $\kappa$ such that for some
$\sigma - k$--linked p.o. $\PP$, Martin's axiom $MA_\kappa$ fails for $\PP$.

\begin{slem} \label{firstlem}
Let $\PP$ be $\sigma - 2^k$--linked, and assume $\dot\phi$ 
is a $\PP$--name for a function $ \bigcup_i 2^{ik} \to 2^k$.
Then there is a countable set $\Psi$ of functions $ \bigcup_i 2^{ik} \to 2^k$ such that
whenever $g \in\twoom$ is such that for all $\psi\in\Psi$ there are infinitely many $i$
with $\psi (g\re ik) = g \re [ik , (i+1) k)$, then
\[ \forces ``\mbox{there are infinitely many } i  \mbox{ with } \dot\phi (g \re ik) = g \re [ik, (i+1) k)." \]
\end{slem}

\begin{proof}
Assume $\PP = \bigcup_n P_n$ where each $P_n$ is $2^k$--linked.
Define $\psi_n  : \bigcup_i 2^{ik} \to 2^k$ such that, for each $\sigma \in 2^{ik}$,
$\psi_n (\sigma)$ is a $\tau$ such that no $p \in P_n$ forces $\dot \phi (\sigma) \neq \tau$.
(Such a $\tau$ clearly exists. For otherwise, for each $\tau \in 2^k$ we could find $p_\tau \in P_n$
forcing $\dot\phi (\sigma) \neq\tau$. Since $P_n$ is $2^k$--linked, the $p_\tau$ would have a
common extension which would force $\dot\phi (\sigma) \notin 2^k$, a contradiction.)
Let $\Psi = \{ \psi_n ; \; n\in\omega\}$.

Now choose $g \in\twoom$ such that for all $\psi\in\Psi$ there are infinitely many $i$
with $\psi (g\re ik) = g \re [ik , (i+1) k)$. 
Fix $i_0$ and $p \in \PP$. There is $n$ such that $p \in P_n$. 
We can find $i \geq i_0$ such that $\psi_n (g\re ik) = g \re [ik , (i+1) k)$.
By definition of $\psi_n$, there is $q \leq p$ such that $q \forces \dot \phi (g\re ik) = \psi_n (g\re ik)$.
Thus $q \forces \dot \phi (g\re ik) = g \re [ik , (i+1) k)$, as required.
\end{proof}

\begin{slem}   \label{secondlem}
Let $\la \PP_n , \dot \QQ_n ; \; n\in\omega \ra$ be a finite support iteration, and assume $\dot\phi$ 
is a $\PP_\omega$--name for a function $ \bigcup_i 2^{ik} \to 2^k$.
Also assume for each $n$ and each $\PP_n$--name $\dot \phi_n$ for a function $ \bigcup_i 2^{ik} \to 2^k$,
there is a countable set $\Psi_n$ of functions $ \bigcup_i 2^{ik} \to 2^k$ such that
$\forall g \in\twoom$, if $\forall\psi\in\Psi_n \; \exists^\infty i \; (\psi (g\re ik) = g \re [ik , (i+1) k))$, then
\[ \forces_n ``\exists^\infty i  \; ( \dot\phi_n (g \re ik) = g \re [ik, (i+1) k))." \]
Then there is a countable set $\Psi$ of functions $ \bigcup_i 2^{ik} \to 2^k$ such that
$\forall g \in\twoom$, if $\forall\psi\in\Psi \; \exists^\infty i \; (\psi (g\re ik) = g \re [ik , (i+1) k))$, then
\[ \forces_\omega ``\exists^\infty i  \; ( \dot\phi (g \re ik) = g \re [ik, (i+1) k))." \]
\end{slem}

\begin{proof}
This is a standard argument which we leave to the reader.
\end{proof}

\begin{slem}   \label{thirdlem}
Let $\PP$ be a p.o. of size $\kappa$, and assume $\dot\phi$ 
is a $\PP$--name for a function $ \bigcup_i 2^{ik} \to 2^k$.
Then there is a set $\Psi$ of size $\kappa$ of functions $ \bigcup_i 2^{ik} \to 2^k$ such that
$\forall g \in\twoom$, if $\forall\psi\in\Psi \; \exists^\infty i \; (\psi (g\re ik) = g \re [ik , (i+1) k))$, then
\[ \forces_\omega ``\exists^\infty i  \; ( \dot\phi (g \re ik) = g \re [ik, (i+1) k))." \]
\end{slem}

\begin{proof}
This is well--known and trivial.
\end{proof}

\medskip

Using the first two of these three lemmata we see that if we iterate $\sigma - 2^k$--linked forcing over a model $V$
containing a family $\F \sub \twoom$ such that
\begin{itemize}
\item[$(\star)$] for all countable sets $\Psi$ of functions $\bigcup_i 2^{ik} \to 2^k$ there is $g \in \F$
   with $\exists^\infty i \; (\psi (g\re ik) = g \re [ik , (i+1) k))$, 
\end{itemize}
then $\F$ still satisfies $(\star)$ in the final extension.
We also have

\begin{slem}   \label{forthlem}
If $\F$ satisfies $(\star)$, then $\ee_2^\const (k) \leq |\F|$.
\end{slem}

\begin{proof}
Simply note $\F$ is a witness for $\ee_2^\const (k)$. For given a predictor $\pi : \twolom \to 2$,
define $\phi : \bigcup_i 2^{ik} \to 2^k$ by $\phi (\sigma) =$ the unique $\tau \in 2^k$
such that $\pi$ predicts $\sigma\ha \tau$ incorrectly on the whole interval $[ik,
(i+1) k)$ where $|\sigma| = ik$. If $g \in \F$ is such that 
$\exists^\infty i \; (\phi (g\re ik) = g \re [ik , (i+1) k))$, then $\pi$ does not $k$--constantly
predict $g$.
\end{proof}

\medskip

Let $2\leq k$.
The partial order $\PP^k$ for adjoining a generic predictor $k$--constantly predicting all ground model
reals is defined as follows. 
Conditions are triples $(\ell,\sigma,F)$ such that $\ell\in\omega$, $\sigma : \twolom \to 2$ is a finite
partial function, and $F \sub \twoom$ is finite, and such that the following requirements are met:
\begin{itemize}
\item $\dom (\sigma ) = 2^{\leq  \ell}$,
\item $f\re \ell \neq g\re \ell$ for all $f \neq g$ belonging to $F$,
\item $\sigma (f \re \ell) = f(\ell)$ for all $f\in F$.
\end{itemize}
The order is given by: $(m, \tau, G) \leq (\ell, \sigma , F)$ if and only if $m \geq\ell$,
$\tau \supseteq \sigma$, $G \supseteq F$, and for all $f\in F$ and all intervals $I \sub (\ell ,m)$
of length $k$ there is $i\in I$ with $\tau (f\re i) = f(i)$.
This is a variation of a p.o. originally introduced in~\cite{Brta}. 
It has been considered as well by Kada~\cite{Kadun}, who also obtained the following lemma.

\begin{slem}
$\PP^k$ is $\sigma - (2^k -1)$--linked.
\end{slem}

\begin{proof}
Simply adapt the argument from~\cite[Lemma 3.2]{Brta}, or see~\cite[Proposition 3.3]{Kadun}.
\end{proof}

\begin{scor} {\rm (Kada~\cite[Corollary 3.5]{Kadun})}   \label{sec3cor}
$\mm (\sigma - (2^k - 1) - \linked) \leq \ee_2^\const (k)$.
\end{scor}

We are ready to prove a result which is dual to Theorem~\ref{Sacksprod}.

\begin{sthm}  \label{3mainthm}
Let $\la \kappa_k ; \; 2 \leq k \in \omega \ra$ be a sequence of uncountable regular cardinals with
$\kappa_k \leq \kappa_{k+1}$. Also assume $\lambda = \lambda^{< \lambda}$ is above the $\kappa_k$.
Then there is a generic extension satisfying $\ee_2^\const (k) = \kappa_k$ for all $k$ and
$\cc = \lambda$. We may also get $\mm (\sigma - (2^k - 1) - \linked) = \kappa_k$ for all $k$.
\end{sthm}

\begin{proof}
Let $\la \PP_\alpha, \dot \QQ_\alpha ; \; \alpha < \lambda \ra$ be a finite support iteration of ccc
forcing such that each factor $\dot \QQ_\alpha$ is forced to be a $\sigma - (2^k -1)$--linked
forcing notion of size less than $\kappa_k$ for some $k \geq 2$. Also guarantee we take care
of all such forcing notions by a book--keeping argument.
Then $\mm (\sigma - (2^k - 1) - \linked) \geq \kappa_k$ is straightforward.
In view of Corollary~\ref{sec3cor} it suffices to prove $\ee_2^\const (k) \leq \kappa_k$ for all $k$.
So fix $k$. Note that in stage $\kappa_k$ of the iteration we adjoined a family $\F$
of size $\kappa_k$ satisfying $(\star)$ above with {\it countable} 
replaced by {\it less than $\kappa_k$}. Show by induction on the remainder of the iteration that
$\F$ continues to satisfy this version of $(\star)$. The limit step is taken care of by Lemma~\ref{secondlem}.
For the successor step, in case $\dot\QQ_\alpha$ is $\sigma - 2^\ell $--linked for some $\ell \geq k$,
use Lemma~\ref{firstlem}, and in case it is not $\sigma - 2^k$--linked
(and thus of size less than $\kappa_k$), use Lemma~\ref{thirdlem}.
By Lemma~\ref{forthlem}, $\ee_2^\const (k) \leq \kappa_k$ follows.
\end{proof}

\medskip

By somewhat changing the above proof, we can dualize Kamo's $CON (\vv^\const_2 > \Cof N)$
(and thus answer a question of his, see~\cite{Kata}), and reprove his result as well.

\begin{sthm}  \label{Kamoques}
\begin{alphenumerate}
\item $\ee^\const_2 < \Add N$ is consistent; in fact, given $\kappa < \lambda = \lambda^{<\kappa}$
regular uncountable, there is a p.o. $\PP$ forcing $\ee^\const_2 = \kappa$ and $\Add N = \cc = \lambda$.
\item {\rm (Kamo, \cite{Ka00})} $\vv^\const_2 > \Cof N$ is consistent; in fact, given $\kappa$ regular
uncountable and $\lambda = \lambda^\omega > \kappa$, there is a p.o. $\PP$ forcing $\vv_2^\const = \cc
=\lambda$ and $\Cof N = \kappa$.
\end{alphenumerate}
\end{sthm}

\begin{proof}
(a) Let $\la \PP_\alpha, \dot \QQ_\alpha ; \; \alpha < \lambda \ra$ be a finite support iteration of ccc
forcing such that 
\begin{itemize}
\item for even $\alpha$, $\forces_\alpha \dot \QQ_\alpha$ is amoeba forcing,
\item for odd $\alpha$, $\forces_\alpha \dot\QQ_\alpha$ is a subforcing of some $\PP^k$ of size less than $\kappa$.
\end{itemize}
Guarantee that we go through all such subforcings by a book--keeping argument.
Then $\ee_2^\const \geq \kappa$ is straightforward, as is $\Add N = \cc = \lambda$.
Now note that amoeba forcing is $\sigma - \infty$--linked (like random forcing).
Therefore we can apply Lemmata~\ref{firstlem}, \ref{secondlem}, and~\ref{thirdlem} for all $k$
simultaneously, and see that there is a family $\F$ of size $\kappa$ 
which satisfies the appropriate modified version of
$(\star)$ (such a family is adjoined after the first $\kappa$ stages of the iteration).

(b) First add $\lambda$ many Cohen reals. Then make a $\kappa$--stage finite support iteration of amoeba
forcing. Again, $\Cof N = \kappa$ is clear. $\vv^\const_2 = \cc = \lambda$ follows from
Lemmata~\ref{firstlem} and~\ref{secondlem} using standard arguments.
\end{proof}

\medskip

One can even strengthen Theorem~\ref{3mainthm} in the following way.
Say a p.o. $\PP$ {\it satisfies property $K_k$}
if for all uncountable $X \sub \PP$ there is $Y \sub X$ uncountable such that any $k$
many elements from $Y$ have a common extension.
Property $K_k$ is a weaker relative of $\sigma -k$--linkedness.
Let $\mm (K_k)$ denote the least cardinal $\kappa$ such that $MA_\kappa$ fails for property $K_k$
p.o.'s.

\begin{slem}
Assume $CH$. $\PP^k$ does not have property $K_{2^k}$. In fact no property $K_{2^k}$ p.o.
adds a predictor which $k$--constantly predicts all ground model reals.
\end{slem}

\begin{proof}
List all predictors as $\{ \pi_\alpha ; \; \alpha < \omega_1 \}$. Choose reals
$f_\alpha \in \twoom$ such that $\pi_\alpha$ does not $k$--constantly predict $f_\beta$
for $\beta \geq \alpha$. Let $X = \{ f_\alpha ; \; \alpha < \omega_1 \}$.

Let $\PP$ be property $K_{2^k}$. Also let $\dot\pi$ be a $\PP$--name for a predictor.
Assume there are conditions $p_\alpha \in \PP$ such that
$p_\alpha \forces ``\dot\pi\;\;\; k$--constantly predicts $f_\alpha$ from $m_\alpha$ onwards."
Without loss $m_\alpha = m$ for all $\alpha$, and any $2^k$ many $p_\alpha$
have a common extension.
Let $T \sub \twolom$ be the tree of initial segments of
members of $X$. Given $\sigma \in T$ with $|\sigma| \geq m$, let $A_\sigma^k = \{ \tau \in T ; \; \sigma \subset \tau$
and $|\tau| = |\sigma| + k \}$. Note that if $|A_\sigma^k | < 2^k$ for all such $\sigma$,
then we could construct a predictor $\pi$ $k$--constantly predicting all of $X$ past $m$ as in
the proof of Theorem~\ref{sec1mainthm}. So there is $\sigma \in T$ with $|A^k_\sigma | = 2^k$.
Find $\alpha_0 , ... , \alpha_{2^k-1}$  such that $A^k_\sigma = \{ f_{\alpha_i} \re |\sigma| + k ; \; 
i < 2^k \}$ and notice
that a common extension of the $p_{\alpha_i}$ forces a contradiction.
\end{proof}

\medskip

Note that some assumption is necessary for the above result for $MA_{\aleph_1}$ implies all
p.o.'s have property $K_k$ for all $k$. We now get

\begin{sthm}
Assume $CH$. Let $2\leq k < \omega$. Then there is a generic extension satisfying
$\ee_2^\const (k) = \aleph_1$ and $\mm (K_{2^k}) = \aleph_2$.
\end{sthm}

\begin{proof}
Use the lemma and the folklore fact that the iteration of property $K_\ell$ p.o.'s is 
property $K_\ell$.
\end{proof}

\medskip

Since we saw in Corollary~\ref{sec3cor} that  $\ee_2^\const (k) \geq \mm (\sigma - (2^k - 1)
- \linked)$. one may ask, on the other hand, whether $\ee_2^\const (k) > \mm (\sigma - (2^k - 1)
- \linked)$ is consistent. This, however, is easy, for the forcing $\PP^k$ is Suslin ccc~\cite{BJ95}
while it is well--known that iterating Suslin ccc forcing keeps numbers like $\mm (\sigma - (2^k - 1)
- \linked)$ small (it even keeps the splitting number $\sss$ small).

We close this section with a few questions. We have no dual result for Theorem~\ref{2mainthm} so far.

\begin{sques}
Is $\ee_2^\const > \sup \{ \ee_2^\const (k) ; \; k < \omega \}$ consistent?
\end{sques}

\begin{sques}
Can $\ee_2^\const$ have countable cofinality?
\end{sques}

By Theorem~\ref{3mainthm}, either of these two questions must have a positive answer.
In fact, in view of the proof of Theorem~\ref{Kamoques}, $\ee_2^\const$ must
be \begin{itemize}
\item either $\max \{ \kappa_k ; \; k\in\omega \}$ (in case the set has a max),  \item or 
$\sup \{ \kappa_k ; \; k\in\omega \}$ or its successor (in case the set has no max)
\end{itemize}
in the model of Theorem~\ref{3mainthm}.


\end{document}